\documentclass[letterpaper, 10pt, conference]{IEEEconf}
\IEEEoverridecommandlockouts 

\usepackage{amsmath,amssymb}
\usepackage{epsfig,subfigure,url,warmread}
\usepackage[all,import]{xy}

\newcommand{\norm}[1]{\ensuremath{\left\| #1 \right\|}}
\newcommand{\bracket}[1]{\ensuremath{\left[ #1 \right]}}
\newcommand{\braces}[1]{\ensuremath{\left\{ #1 \right\}}}
\newcommand{\parenth}[1]{\ensuremath{\left( #1 \right)}}
\newcommand{\refeqn}[1]{(\ref{eqn:#1})}
\newcommand{\reffig}[1]{Fig. \ref{fig:#1}}
\newcommand{\tr}[1]{\mbox{tr}\ensuremath{\negthickspace\bracket{#1}}}
\newcommand{\deriv}[2]{\ensuremath{\frac{\partial #1}{\partial #2}}}
\newcommand{\SO}{\ensuremath{\mathrm{SO(3)}}}
\newcommand{\T}{\ensuremath{\mathrm{T}}}
\newcommand{\so}{\ensuremath{\mathfrak{so}(3)}}

\renewcommand{\Re}{\ensuremath{\mathbb{R}}}
\renewcommand{\S}{\ensuremath{\mathbb{S}}}

\renewcommand{\theenumi}{\roman{enumi}}
\renewcommand{\labelenumi}{(\theenumi)}

\title{\LARGE \bf
Attitude Maneuvers of a Rigid Spacecraft in a Circular Orbit}

\author{ \parbox{3 in}{\centering Taeyoung Lee\authorrefmark{1}\authorrefmark{2}, N. Harris McClamroch\authorrefmark{2}\\
         Department of Aerospace Engineering\\
         University of Michigan, Ann Arbor, MI 48109\\
         {\tt\small \{tylee, nhm\}@umich.edu}}
         \hspace*{ 0.5 in}
         \parbox{3 in}{\centering Melvin Leok\authorrefmark{1}\\
         Department of Mathematics\\
        University of Michigan, Ann Arbor, MI 48109\\
         {\tt\small mleok@umich.edu}}
        \thanks{\textsuperscript{\footnotesize\ensuremath{*}}This research has been supported in part by NSF under grant DMS-0504747, and by a grant from the Rackham Graduate School, University of Michigan.}
        \thanks{\textsuperscript{\footnotesize\ensuremath{\dagger}}This research has been supported in part by NSF under grant ECS-0244977.}
}

\begin{document}
\allowdisplaybreaks
\maketitle \thispagestyle{empty} \pagestyle{empty}

\begin{abstract}
A global model is presented that can be used to study attitude
maneuvers of a rigid spacecraft in a circular orbit about a large
central body.   The model includes gravity gradient effects that
arise from the non-uniform gravity field and characterizes the
spacecraft attitude with respect to the uniformly rotating local
vertical local horizontal coordinate frame. An accurate
computational approach for solving a nonlinear boundary value
problem is proposed, assuming that control torque impulses can be
applied at initiation and at termination of the maneuver. If the
terminal attitude condition is relaxed, then an accurate
computational approach for solving the minimal impulse optimal
control problem is presented. Since the attitude is represented by a
rotation matrix, this approach avoids any singularity or ambiguity
arising from other attitude representations such as Euler angles or
quaternions.
\end{abstract}

\section{Introduction}
The attitude dynamics of an uncontrolled rigid spacecraft in a
circular orbit about a large central body, including gravity
gradient effects, have been extensively studied; see
\cite{bk:hughes,bk:wie}.   There are 24 distinct relative equilibria
for which the principal axes are exactly aligned with the local
vertical local horizontal (LVLH) axes, and the spacecraft angular
velocity is identical to the orbital angular velocity of the LVLH
coordinate frame. Linear rotational equations of motion that
describe small perturbations from any relative equilibrium solutions
are well known.   Linear attitude control of a rigid spacecraft in a
circular orbit, including linear gravity gradient effects, has also
been addressed in \cite{bk:wie}. However, linear controllers have the
limitation that they are only applicable to small attitude change
maneuvers.

The emphasis in this paper is on large angle attitude maneuvers of a
rigid spacecraft in a circular orbit about a large central body,
including gravity gradient effects.   A nonlinear, globally defined
model is introduced. This model describes the attitude of the
spacecraft, relative to the uniformly rotating  LVLH coordinate
frame, by a rotation matrix.   The model includes gravity gradient
 terms that reflect the rotation of the LVLH frame, and terms
that reflect control input torques.   In the problems studied in
this paper, independent impulsive control torques can be applied
about each principal axis. Gravity gradient moments are significant
in Earth orbits for orbit altitudes between $400 \,\mathrm{km}$ and
$40,000\, \mathrm{km}$ and for attitude maneuver times that are not
small compared with the orbital period.

We study open loop attitude maneuvers that can be accomplished by
using two impulsive torque controls, one occurring at the initial
time and one occurring at the terminal time of the maneuver.   The
attitude motion in between the initial time and the final time is
uncontrolled.    Two classes of attitude maneuver problems are
studied.   In Section \ref{sec:tpbvp}, a terminal attitude is
specified so that there is a unique impulse sequence satisfying the
boundary conditions.  This problem can be solved computationally use
a root finding algorithm.   In Section \ref{sec:opt}, the specified
terminal attitude condition is relaxed such that there are many
impulse sequences that satisfy the boundary conditions.  In this
case we seek the minimum total impulse sequence that satisfies the
boundary conditions.   This problem can be solved computationally
using a constrained minimization algorithm.

For both classes of attitude maneuver problems, our contribution is to
demonstrate how sensitivity derivatives, used in the computational
algorithms, can be determined effectively and accurately such that
they satisfy the global geometry of the problem. Since the attitude
is represented by a rotation matrix in the special orthogonal group
\SO, and the sensitivity derivatives are expressed in terms of the
Lie algebra \so, this approach completely avoids singularities or
ambiguities that arise from other representations of the rotation group, such as Euler angles
and quaternions.

\section{Rigid spacecraft models in a circular orbit}

\renewcommand{\xyWARMinclude}[1]{\includegraphics[width=0.45\textwidth]{#1}}
\begin{figure}[t]
{\vspace*{-0.5cm}\scriptsize\selectfont
$$\begin{xy}
\xyWARMprocessEPS{coordinate}{eps}
\xyMarkedImport{}
\xyMarkedMathPoints{1-9}
\end{xy}$$\vspace*{-1cm}}
\caption{Coordinate frames}\label{fig:coord}\vspace*{-0.5cm}
\end{figure}
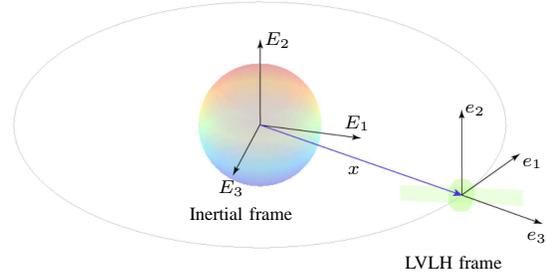

We assume that a rigid spacecraft is on a circular orbit with a
constant orbital angular velocity $\omega_0\in\Re$.
In this section, the continuous equations of motion and a geometric numerical integrator, referred
to as Lie group variational integrator, for the attitude maneuver of
a spacecraft in a circular orbit are given following the results of
\cite{pro:CCA05}. We identify $\T\SO\simeq\SO\times\so$ by left
translation, and we identify $\so\simeq\Re^3$ by an isomorphism
$S(\cdot):\Re^3\mapsto \so$.

\subsubsection*{Continuous equations of motion} We define three
rotation matrices in $\SO$;
\begin{align*}
R^{bi} &: \text{ from the body fixed frame to the
inertial frame},\\
R^{li} &: \text{ from the LVLH frame to the inertial frame},\\
R^{bl} &: \text{ from the body fixed frame to the LVLH frame},
\end{align*}
where the inertial frame and the LVLH frame are illustrated in
\reffig{coord}. Thus, $R^{bl}=R^{{li}^T}R^{bi}$.

The on-orbit spacecraft equations of motion are given by
\begin{gather}
    \dot{\Pi}+\Omega\times \Pi = M^g, \label{eqn:dotPi}\\
    \dot{R}^{bi}=R^{bi}S(\Omega),\label{eqn:dotRbi}\\
    \dot{R}^{li}=R^{li}S(\omega_0 e_2),\label{eqn:dotRli}\\
    \dot{R}^{bl}=R^{bl}S(\Omega-\omega_0 {R^{{bl}^T}}e_2)\label{eqn:dotRbl},
\end{gather}
where $\Pi,\Omega\in\Re^3$ are the angular momentum and the angular
velocity of the spacecraft expressed in the body fixed frame,
respectively, and $M^g\in\Re^3$ is the gravity gradient moment. The
isomorphism between $\Re^3$ and $\so$ is defined such that
$S(x)y=x\times y$ for any $x,y\in\Re^3$.
Since the orbital angular velocity $\omega_0$ is constant, the
solution of \refeqn{dotRli} is given by $R^{li}(t)=R^{li}(0)
e^{S(\omega_0 e_2)t}$.

\subsubsection*{Gravity gradient moment} The gravity gradient moment
is derived in \cite{bk:wie}. We present an alternative way to obtain
the gravity gradient moment directly using the gravity potential;
\begin{align*}
U=-\int_{\mathcal{B}} \frac{GM}{\norm{x+R^{bi}\rho}}\; dm,
\end{align*}
where $x\in\Re^3$ is the position of the spacecraft in the inertial
frame, and $\rho\in\Re^3$ is a vector from the center of mass of the
spacecraft to a mass element in the body fixed frame. $G$ is the
gravitational constant and $M$ is the mass of the Earth.

>From \cite{pro:CCA05}, the gravity gradient moment $M^g$ can be
determined by using the following relationship;
\begin{align}
S(M^g)=\deriv{U}{R^{bi}}^T R^{bi}
-R^{{bi}^T}\deriv{U}{R^{bi}}.\label{eqn:SMg}
\end{align}
We derive a closed form for $M^g$ from \refeqn{SMg}, by assuming
that the spacecraft is on a circular orbit so that the norm of $x$
is constant, and the size of the spacecraft is much smaller than the
size of the orbit.

As shown in \reffig{coord}, the coordinate of the spacecraft
position in the LVLH frame is
 $r_0e_3$, where $r_0\in\Re$ is the radius of the
circular orbit, and $e_3=\bracket{0,\,0,\,1}^T$. Therefore, the
position of the spacecraft in the inertial frame is given by $x=r_0
R^{li}e_3$. Using this expression,
\begin{align*}
\deriv{U}{R^{bi}} & = \int_{\mathcal{B}} \frac{GM\, r_0 R^{li}\,e_3
\rho^T}{\norm{\,r_0 e_3+R^{bl}\rho\,}^3}\;dm,\\
& = \frac{GM}{r_0}\int_{\mathcal{B}} 
\frac{\parenth{R^{li}e_3
\hat\rho^T}\frac{\norm{\rho}}{r_0}}{\bracket{1+2\parenth{e_3^T
R^{bl}\hat{\rho}}\frac{\norm{\rho}}{r_0}+\frac{\norm{\rho}^2}{r_0^2}}^{\frac{3}{2}}}\;dm,
\end{align*}
where $\hat\rho=\frac{\rho}{\norm{\rho}}\in\Re^3$ is the unit vector
along the direction of $\rho$. Assume that the size of the
spacecraft is significantly smaller than the size of the orbit, i.e.
$\frac{\norm{\rho}}{r_0}\ll 1$. Using a Taylor series expansion, we
obtain the 2nd order approximation.
\begin{align*}
\deriv{U}{R^{bi}} & = \frac{GM}{r_0}\int_{\mathcal{B}} 
R^{li}e_3 \hat\rho^T\braces{\frac{\norm{\rho}}{r_0}
-3 e_3^T R^{bl}\hat{\rho} \frac{\norm{\rho}^2}{r_0^2}}\;dm.
\end{align*}
Since the body fixed frame is located at the mass center of the
spacecraft, $\int_{\mathcal{B}}\rho\;dm=0$. Therefore, the first
term in the above equation vanishes. Because $e_3^T R^{bl}\rho$ is a
scalar quantity, we can rewrite the above equation as
\begin{align}
\deriv{U}{R^{bi}}
& = -3 \omega_0^2R^{li}e_3e_3^T
R^{bl}\parenth{\frac{1}{2}\tr{J}I_{3\times 3}-J},\label{eqn:UR}
\end{align}
where $\omega_0=\sqrt{\frac{GM}{r_0^3}}\in\Re$ is the orbital
angular velocity, and $J\in\Re^{3\times 3}$ is the moment of inertia
matrix of the spacecraft. Substituting \refeqn{UR} into
\refeqn{SMg}, and
using the property $S(x\times y)=yx^T-xy^T$ for $x,y\in\Re^3$, we
obtain an expression for the gravity gradient moment as follows.
\begin{align}
    M^g & = 3\omega_0^2 R^{{bl}^T}e_3 \times J R^{{bl}^T}
    e_3.\label{eqn:Mg}
\end{align}

\subsubsection*{Discrete equations of motion} In the continuous
equations of motion, the structure of \refeqn{dotRbi},
\refeqn{dotRli}, and \refeqn{dotRbl} ensures that $R^{bi}$,
$R^{li}$, and $R^{bl}$ evolve on the special orthogonal group,
$\SO$. However, general numerical integration methods, including the
popular Runge-Kutta methods, do not preserve the orthogonality
property of this group. For example, if we integrate \refeqn{dotRbl}
by a typical Runge-Kutta scheme, the quantity $R^{{bl}^T}R^{bl}$
inevitably drifts from the identity matrix as the simulation time
increases.

The rotation matrix is commonly parameterized by Euler angles or
quaternions. These attitude kinematics equations can be numerically
integrated and are used to recompute the rotation matrix. However,
Euler angles are not global expressions of the attitude since they
have associated singularities. The analytical expressions for
sensitivities are hard to develop since many trigonometric terms are
encountered. Quaternions have no singularity, but quaternions must
lie on the three sphere $\S^3$. General numerical integration
methods do not preserve the unit length of a quaternion. Therefore,
quaternions have the same numerical drift problem as rotation
matrices. Furthermore, quaternions, which are diffeomorphic to
$\mathrm{SU(2)}$, double covers $\SO$. So there are inevitable
ambiguities in expressing the attitude.

These cause significant inaccuracies in numerical simulations based on quaternions and Euler angles. In particular, the gravity gradient moment, as given in
\refeqn{Mg}, depends on $R^{bl}$ directly, and consequently, errors
in computing $R^{bl}$ cause errors in the gravity gradient moment.
These effects are more pronounced when the simulation time is large.

Lie group variational integrators preserve the orthonormal structure
of $\SO$ without any reprojection or parametrization. They also
conserve the momentum map, and the symplectic property of rigid body
dynamics. In addition, the total energy is well-behaved, as it only
oscillates in a bounded fashion about its true value. So, Lie group
variational integrators are geometrically exact. Using the results
given in \cite{pro:CCA05}, a Lie group variational integrator for
the attitude dynamics of a spacecraft in a circular orbit are given
by
\begin{gather}
\Pi_{k+1} = F_k^T \Pi_k +\frac{h}{2} F_k^T M_k^g +\frac{h}{2}
M_{k+1}^g,\label{eqn:updatePi}\\
h S(\Pi_k+\frac{h}{2} M_k^g) = F_k J_d - J_d
F_k^T,\label{eqn:findf}\\
R_{k+1}^{bi} = R_k^{bi} F_k,\label{eqn:updateRbi}\\
R_{k+1}^{li}=R_k^{li} e^{S(\omega_0 e_2 )h},\label{eqn:updateRli}\\
R_{k+1}^{bl}= e^{-S(\omega_0 e_2) h}R_k^{bl}
F_k,\label{eqn:updateRbl}
\end{gather}
where the subscript $k$ denotes variables at the $k$th time step,
and $h\in\Re$ is the integration step size. The matrix
$J_d\in\Re^{3\times 3}$ is a nonstandard moment of inertia matrix
defined by $J_d=\frac{1}{2}\tr{J}I_{3\times 3}-J$.

The matrix $F_k=R_k^{bi^T}R_{k+1}^{bi}\in\SO$ is the relative
attitude between integration steps, and it is obtained by solving
\refeqn{findf}.
Since $F_k$ and $e^{S(\omega_0 e_2)h}$ are in $\SO$, and $\SO$ is closed under matrix multiplication, $R^{bi}_k$, $R^{li}_k$,
and $R^{bl}_k$ evolve in $\SO$ automatically for all $k$ according
to \refeqn{updateRbi}, \refeqn{updateRli}, and \refeqn{updateRbl}.
The actual computation of $F_k$ is done in the Lie algebra $\so$ of
dimension 3, and the rotation matrices are updated by
multiplication with the exponential of a skew-symmetric matrix. So, Lie group variational integrators are numerically
efficient, and there is no excessive computational burden in updating the
9 elements of the rotation matrix. The properties of these discrete
equations of motion are discussed more explicitly in
\cite{pro:CCA05} and \cite{jo:CMAME05}.

Since $\omega_0$ is constant, the solution of \refeqn{updateRli} is
given by
\begin{align}
R_{k}^{li}=R_0^{li} e^{S(\omega_0 e_2 )kh}.\label{eqn:Rlik}
\end{align}

\section{Spacecraft attitude maneuvers}\label{sec:tpbvp}

\subsection{Problem formulation}
A two point boundary value problem is formulated for spacecraft
rest-to-rest maneuver between two given attitudes for a fixed
maneuver time. This is a Lambert type boundary value problem on
$\SO$. The initial attitude and the desired terminal attitude are
expressed as rotation matrices with respect to the LVLH frame,
namely $R^{bl}_{0},\,R^{bl}_{N_d}\in\SO$. Two impulsive control
torques are applied at the initial time and the terminal time. We
assume that the control torques are purely impulsive, which means
that each impulse changes the angular velocity of the spacecraft
instantaneously, but it does not have any effect on the attitude of
the spacecraft at that instant. The rotational motion of the
spacecraft between the initial time and the terminal time is
uncontrolled. This assumption is reasonable when the operating time
of the spacecraft control moment actuator is much smaller than the
fixed maneuver time.

We define this boundary value problem directly on $\SO$ instead of
using parameterizations such as Euler angles and quaternions. This
approach allows us to define the attitude of the spacecraft globally
 without singularities and ambiguities. We use the discrete equations of motion given in
\refeqn{updatePi}--\refeqn{updateRbl} for the problem formulation
and for the following analysis.

We transform this two point boundary value problem into a nonlinear
root finding problem. For a given initial angular momentum $\Pi_0$
and an initial attitude of the spacecraft $R^{bl}_0$, the terminal
angular momentum $\Pi_N$ and the terminal attitude $R^{bl}_N$ are
determined by the discrete equations of motion. By choosing the initial
angular momentum so that the terminal attitude of the
spacecraft is equal to the desired attitude, i.e.
$R^{bl}_N=R^{bl}_{N_d}$, we obtain the initial impulse
 and the terminal impulse. Thus, the nonlinear boundary value problem for the spacecraft
attitude maneuver is formulated as
\begin{gather*}
\text{given}  : R^{bl}_0, R^{bl}_{N_d}, N\\
\text{find}  : \Pi_0 \\
\text{such that } R^{bl}_{N}=R^{bl}_{N_d}  \text{ subject to
\refeqn{updatePi}--\refeqn{updateRbl},}
\end{gather*}
where $N\in\mathbb{N}$ is the number of integration steps determined
by $N=\frac{T}{h}$ for the fixed maneuver time $T$ and the fixed
integration step size $h$.

\subsection{Computational approach}
We solve a sequence of linear boundary value problems whose
solutions converge to the solution of the nonlinear boundary value
problem.

\subsubsection*{Linearization} The equations of motion are linearized about a given trajectory, and they are expressed in terms of the Lie algebra $\so$.
Consider small perturbations from a given trajectory denoted by
$\Pi_k^\epsilon,\,R^{bi,\epsilon}_k,R^{bl,\epsilon}_k,F_k^\epsilon$:
\begin{align}
\Pi_k^\epsilon & = \Pi_k + \epsilon\delta\Pi_k,\label{eqn:Pie}\\
R^{bi,\epsilon}_k &= R^{bi}_k + \epsilon \delta R^{bi}_k +
\mathcal{O}(\epsilon^2),\label{eqn:Rbie}\\
R^{bl,\epsilon}_k &= R^{bl}_k + \epsilon \delta R^{bl}_k +
\mathcal{O}(\epsilon^2),\label{eqn:Rble}\\
F^{\epsilon}_k &= F_k + \epsilon \delta F_k +
\mathcal{O}(\epsilon^2),\label{eqn:Fe}
\end{align}
where $\epsilon\in\Re$. Since the orbital angular velocity is
constant, $\delta R^{li}_k=0$.

The infinitesimal variation of the angular momentum $\delta \Pi_k$ can be expressed in $\Re^3$. The variation of the rotation matrix
$R^{bi,\epsilon}_k\in\SO$ can be expressed as
\begin{align*}
R^{bi,\epsilon}_k = R^{bi}_k e^{\epsilon S(\zeta_k)},
\end{align*}
where $\zeta_k\in\Re^3$ and $S(\zeta_k)\in\so$ is a skew-symmetric matrix. Since the map $S(\cdot)$ is an isomorphism between $\so$ and $\Re^3$, $\zeta_k$ is well defined. Then, the
infinitesimal variation $\delta R^{bi}_k$ is given by
\begin{align}
\delta R^{bi}_k =
\frac{d}{d\epsilon}\bigg|_{\epsilon=0}R^{bi,\epsilon}_k=R^{bi}_kS(\zeta_k).\label{eqn:delRbik}
\end{align}
$\delta F_k$ is obtained from definition \refeqn{updateRbi}, and
\refeqn{delRbik} as
\begin{align}
\delta F_k&=
\frac{d}{d\epsilon}\bigg|_{\epsilon=0}R_k^{bi,\epsilon^T}R_{k+1}^{bi,\epsilon},\nonumber\\
& = -S(\zeta_k) F_k + F_{k}S(\zeta_{k+1}).\label{eqn:delFk}
\end{align}
Since $\delta R^{li}_k=0$, $\delta R^{bl}_k$ is given by
\begin{align}
\delta R^{bl}_k &= R_k^{li^T} \delta R_k^{bi} =
R^{bl}_kS(\zeta_k).\label{eqn:delRblk}
\end{align}
In summary, equations \refeqn{delRbik}, \refeqn{delFk} and
\refeqn{delRblk} describe the variations of rotation matrices in
$\SO$.

Substituting \refeqn{Pie}, \refeqn{Rble}, and \refeqn{Fe} into
\refeqn{updatePi}, and \refeqn{findf}, and ignoring higher-order
terms, the linearized discrete equations of motion are given by
\begin{gather}
\begin{aligned}
\delta \Pi_{k+1} & = \delta F_k^T \Pi_k + F_k^T \delta\Pi_k\\
&\quad +\frac{h}{2} \delta F_k^TM^g_k
+\frac{h}{2}  F_k^T\delta M^g_k
+ \frac{h}{2} \delta M^g_{k+1},
\end{aligned}\label{eqn:updatePid}\\
h S(\delta\Pi_k+\frac{h}{2} \delta M^g_k)
= \delta F_k J_d - J_d \delta F_k^T,\label{eqn:findfd}
\end{gather}
where
\begin{align}
\delta M^g_k & = 3\omega_0^2  \bracket{\delta  R^{bl^T}_ke_3 \times
JR^{bl^T}_ke_3
+ R^{bl^T}_ke_3\times J \delta R^{bl^T}_ke_3}.\label{eqn:delMgk0}
\end{align}

These equations are not in standard form since \refeqn{findfd} is an
implicit equation in $\delta F_k$. By using \refeqn{delRbik},
\refeqn{delFk} and \refeqn{delRblk}, we will obtain an explicit solution of
\refeqn{findfd}, and rewrite the above equations in standard form.

Substituting \refeqn{delRblk}
into \refeqn{delMgk0}, and using the property $S(x)y=-S(y)x$ for all
$x,y\in\Re^3$, $\delta M^g_k$ can be written as
\begin{align}
\delta M^g_k & = 3\omega_0^2 \Big[ -S(JR^{bl^T}_ke_3)S(R^{bl^T}_ke_3)
\nonumber\\& \qquad\qquad +
 S(R^{bl^T}_ke_3) J S(R^{bl^T}_ke_3)\Big]\zeta_k,\nonumber\\
& = \mathcal{M}_k \zeta_k,\label{eqn:delMgk}
\end{align}
where $\mathcal{M}_k\in\Re^{3\times 3}$.

Using \refeqn{delFk} and
the property $S(Rx)=RS(x)R^T$ for $x\in\Re^3$, $R\in\SO$, the right hand side of \refeqn{findfd} is
written as
\begin{align}
\delta  & F_k J_d - J_d \delta F_k^T\nonumber\\
& = -\braces{S(\zeta_k) F_k J_d+J_d F_k^TS(\zeta_k)}\nonumber\\
& \quad +\braces{S(F_k\zeta_{k+1})F_kJ_d + J_dF_k^TS(F_k\zeta_{k+1})}.\label{eqn:RHSfindfd0}
\end{align}

Substituting \refeqn{delMgk} and \refeqn{RHSfindfd0} into \refeqn{findfd}, and using the property
$S(x)A+A^TS(x)=S(\braces{\tr{A}I_{3\times 3}-A}x)$ for $x\in\Re^3$, $A\in\Re^{3\times 3}$, \refeqn{findfd} can
be transformed into an equivalent vector form;
\begin{align*}
h \delta\Pi_k+\frac{h^2}{2} \mathcal{M}_k\zeta_k
=&-\braces{\tr{F_k J_d}I_{3\times 3}-F_k
J_d}\zeta_k\\
&+\braces{\tr{F_k J_d}I_{3\times 3}-F_k J_d}F_k\zeta_{k+1}.
\end{align*}
Multiplying both sides by $F_k^T\braces{\tr{F_k J_d}I_{3\times
3}-F_k J_d}^{-1}$ and rearranging, we obtain
\begin{align}
\zeta_{k+1}&=
F_k^T\bracket{\frac{h^2}{2}\braces{\tr{F_k J_d}I_{3\times 3}-F_k
J_d}^{-1} \mathcal{M}_k+I_{3\times 3}}\zeta_k\nonumber\\
&\quad+hF_k^T\braces{\tr{F_k J_d}I_{3\times 3}-F_k J_d}^{-1}
\delta\Pi_k,\nonumber\\
& \triangleq \mathcal{A}_k \zeta_k+\mathcal{B}_k \delta\Pi_k
,\label{eqn:zetakp}
\end{align}
where $\mathcal{A}_k, \mathcal{B}_k\in\Re^{3\times 3}$.
Substituting \refeqn{zetakp} into \refeqn{delFk}, $\delta
F_k$ is obtained as
\begin{align}
\delta F_k & = -\eta_k F_k + F_k \eta_{k+1},\nonumber\\
& = -S(\zeta_k) F_k + F_k S(\zeta_{k+1}),\nonumber\\
& = -S(\zeta_k) F_k + F_k S(\mathcal{A}_k \zeta_k+\mathcal{B}_k
\delta\Pi_k).\label{eqn:delFks}
\end{align}
Equation \refeqn{delFks} is an explicit solution of \refeqn{findfd}.

Now we rewrite \refeqn{updatePid} in the standard form of linear
discrete equations of motion. Substituting \refeqn{delMgk},
\refeqn{delFks} into \refeqn{updatePid},
and rearranging, we obtain
\begin{align}
\delta \Pi_{k+1} & = \bigg[ F_k^T S( \Pi_k +\frac{h}{2}M^g_k)\braces{-I_{3\times 3}
    + F_k \mathcal{A}_k}\nonumber\\
    & \quad\quad+\frac{h}{2}F_k^T\mathcal{M}_k + \frac{h}{2}\mathcal{M}_{k+1}\mathcal{A}_k\bigg]\zeta_k\nonumber\\
& \!\!\!\!\!\!\! + \bracket{S(F_k^T\braces{\Pi_k +\frac{h}{2}
M^g_k})\mathcal{B}_k
    +F_k^T+\frac{h}{2} \mathcal{M}_{k+1}\mathcal{B}_k}\delta\Pi_k,\nonumber\\
& \triangleq \mathcal{C}_k \zeta_k + \mathcal{D}_k\delta \Pi_k,\label{eqn:delPikp}
\end{align}
where $\mathcal{C}_k,\mathcal{D}_k\in\Re^{3\times 3}$.

In summary, \refeqn{zetakp} and \refeqn{delPikp} are linear discrete
equations equivalent to \refeqn{updatePid}, \refeqn{findfd} and
\refeqn{delMgk0}, and they  can be written as
\begin{align}
\begin{bmatrix}\zeta_{k+1}\\\delta\Pi_{k+1}\end{bmatrix}
& = \begin{bmatrix} \mathcal{A}_k & \mathcal{B}_k \\ \mathcal{C}_k &
\mathcal{D}_k \end{bmatrix}
\begin{bmatrix}\zeta_{k}\\\delta\Pi_{k}\end{bmatrix},\nonumber\\
& \triangleq A_k
\begin{bmatrix}\zeta_{k}\\\delta\Pi_{k}\end{bmatrix},\label{eqn:sl}
\end{align}
where $A_k\in\Re^{6\times 6}$. Equation \refeqn{sl} is the linear
discrete equation for perturbations of the attitude dynamics of a
spacecraft in a circular orbit, expressed in terms of
$\Re^3\simeq\so$. The important feature of \refeqn{sl} is that it is
linearized in such a way that it respects the geometry of the
special orthogonal group $\SO$.

\subsubsection*{Linear boundary value problem}
The solution of \refeqn{sl} is given by
\begin{align}
\begin{bmatrix}\zeta_{N}\\\delta\Pi_{N}\end{bmatrix}& =
\left(\prod_{k=0}^{N-1}A_k\right)
\begin{bmatrix}\zeta_{0}\\\delta\Pi_{0}\end{bmatrix},\nonumber\\
& \triangleq
\begin{bmatrix} \Phi_{11} & \Phi_{12} \\ \Phi_{21} & \Phi_{22}\end{bmatrix}
\begin{bmatrix}\zeta_{0}\\\delta\Pi_{0}\end{bmatrix},\label{eqn:solsl}
\end{align}
where $\Phi_{ij}\in\Re^{3\times 3}$ for $i,j=1,2$.

For the given boundary value problem, $\zeta_0=0$ since the initial
attitude of the spacecraft is given and fixed, and $\delta\Pi_N$ is
free since the terminal angular momentum is compensated by the
terminal impulse. Then, we obtain
\begin{align*}
\zeta_N & = \Phi_{12} \delta\Pi_0.
\end{align*}
This equation provides $\Phi_{12}$, the sensitivity derivative of
the terminal attitude with respect to a change in the initial
angular momentum. It states that for a given trajectory, if we
update the initial angular velocity by $\delta\Pi_0$, then the
terminal attitude is changed from $R^{bl}_N$ to
$R^{bl}_Ne^{S(\zeta_N)}=R^{bl}_Ne^{S(\Phi_{12} \delta\Pi_0)}$.

We choose the change of the initial angular momentum so that the
updated terminal attitude is equal to the desired terminal attitude;
$R^{bl}_Ne^{S(\Phi_{12} \delta\Pi_0)}=R^{bl}_{N_d}$,
or equivalently,
\begin{align}
\delta\Pi_0=\Phi_{12}^{-1}
S^{-1}\!\parenth{\mathrm{logm}\!\parenth{R^{bl^T}_NR^{bl}_{N_d}}},\label{eqn:delPi0}
\end{align}
where $S^{-1}(\cdot):\so\mapsto\Re^3$ is the inverse mapping of
$S(\cdot)$, and $\mathrm{logm}$ denotes the matrix logarithm.
Equation \refeqn{delPi0} provides a solution of the linear boundary
value problem for the attitude dynamics of a spacecraft in a
circular orbit, assuming that $\Phi_{12}$ is invertible.

\subsubsection*{Nonlinear boundary value problem}
The linear boundary value problem is solved successively so that its
solution converges to the solution of the nonlinear boundary value
problem. A numerical algorithm is summarized as follows.

\vspace*{0.1cm}
{
\renewcommand{\theenumi}{\arabic{enumi}}
\renewcommand{\labelenumi}{\theenumi:}
\hrule\vspace*{0.08cm}
\begin{enumerate}
\item Set $\mathrm{Error}=2\epsilon_S$.
\item Guess an initial condition $\Pi_0^{(0)}$.
\item Set $i=0$.
\item \textbf{while} $\mathrm{Error} > \epsilon_S$.
\item Find $\Pi_k^{(i)},R_k^{bl(i)}$ using
$\Pi_0^{(i)}$ and \refeqn{updatePi}, \refeqn{findf},
\refeqn{updateRbl}.
\item Compute the error;\\
$\zeta_N^{(i)}=S^{-1}\!\parenth{\mathrm{logm}\!
\parenth{R_N^{bl(i)^T} R_{N_d}^{bl}}}$,
$\mathrm{Error}=\norm{\zeta_N^{(i)}}$.
\item Update the initial condition; $\Pi_0^{(i+1)}=\Pi_0^{(i)}+c \Phi_{12}^{-1}
\zeta_N^{(i)}$.
\item Set $i=i+1$.
\item \textbf{end while}
\end{enumerate}
\vspace*{0.08cm}
\hrule}
\vspace*{0.1cm}

\noindent Here the superscript $(i)$ denotes the $i$th
iteration, and $\epsilon_S,c\in\Re$ are a stopping criterion and a
scaling factor, respectively.

This computational approach utilizes an exact and efficient method
to compute sensitivity derivatives in the special orthogonal group
$\SO$. The sensitivity derivatives are then used to solve the two
point boundary value problem.

\subsection{Numerical example}
Three spacecraft rest-to-rest maneuvers between relative equilibrium
attitudes are considered. The resulting motions are highly
nonlinear, large angle maneuvers.

The mass, length and time dimensions are normalized. The moment of
inertia of the spacecraft and simulation parameters are chosen as
${J}=\mathrm{diag}\bracket{1,\,2.8,\,2}$, $\epsilon_s = 10^{-14}$,
$c=0.1$, $h=0.001$. Each maneuver is completed in a quarter of the
orbit, ${T}=\frac{\pi}{2}$. The boundary conditions and the
corresponding computed impulsive control are as follows.
\begin{enumerate}
\item Rotational maneuver about the LVLH axis $e_1$:
\begin{gather*}
R_{0}^{bl}=I_{3\times 3},\quad
R_{N_d}^{bl}=\mathrm{diag}\bracket{1, -1, -1},\\
{\Pi}_0=\bracket{2.116,\,1.531,\,-1.782}^T,\\
-{\Pi}_N=\bracket{-2.116,\,1.531,\,1.782}^T.
\end{gather*}
\item Rotational maneuver about the LVLH axes $e_1$ and $e_2$:
\begin{gather*}
R_{0}^{bl}=\mathrm{diag}\bracket{1, -1, -1},\quad
R_{N_d}^{bl}=\begin{bmatrix}
     -1&     0&     0\\
    0&     0&     -1\\
     0&     -1&     0\end{bmatrix},\\
{\Pi}_0=\bracket{-1.323,\,1.798,\,0.932}^T,\\
-{\Pi}_N=\bracket{0.397,\,-1.586,\,-1.310}^T.
\end{gather*}
\item Rotational maneuver about the LVLH axes $e_2$ and $e_3$:
\begin{gather*}
R_{0}^{bl}=\mathrm{diag}\bracket{1, -1, -1},\quad
R_{N_d}^{bl}=\begin{bmatrix}
     0&     1&     0\\
    -1&     0&     0\\
     0&     0&     1\end{bmatrix},\\
{\Pi}_0=\bracket{1.047,\,0.437,\,2.800}^T,\\
-{\Pi}_N=\bracket{-1.416,\,-1.761,\,1.159}^T.
\end{gather*}
\end{enumerate}

\begin{figure}
    \subfigure[Rotation about the LVLH axis $e_1$]{
    \includegraphics[width=0.45\columnwidth]{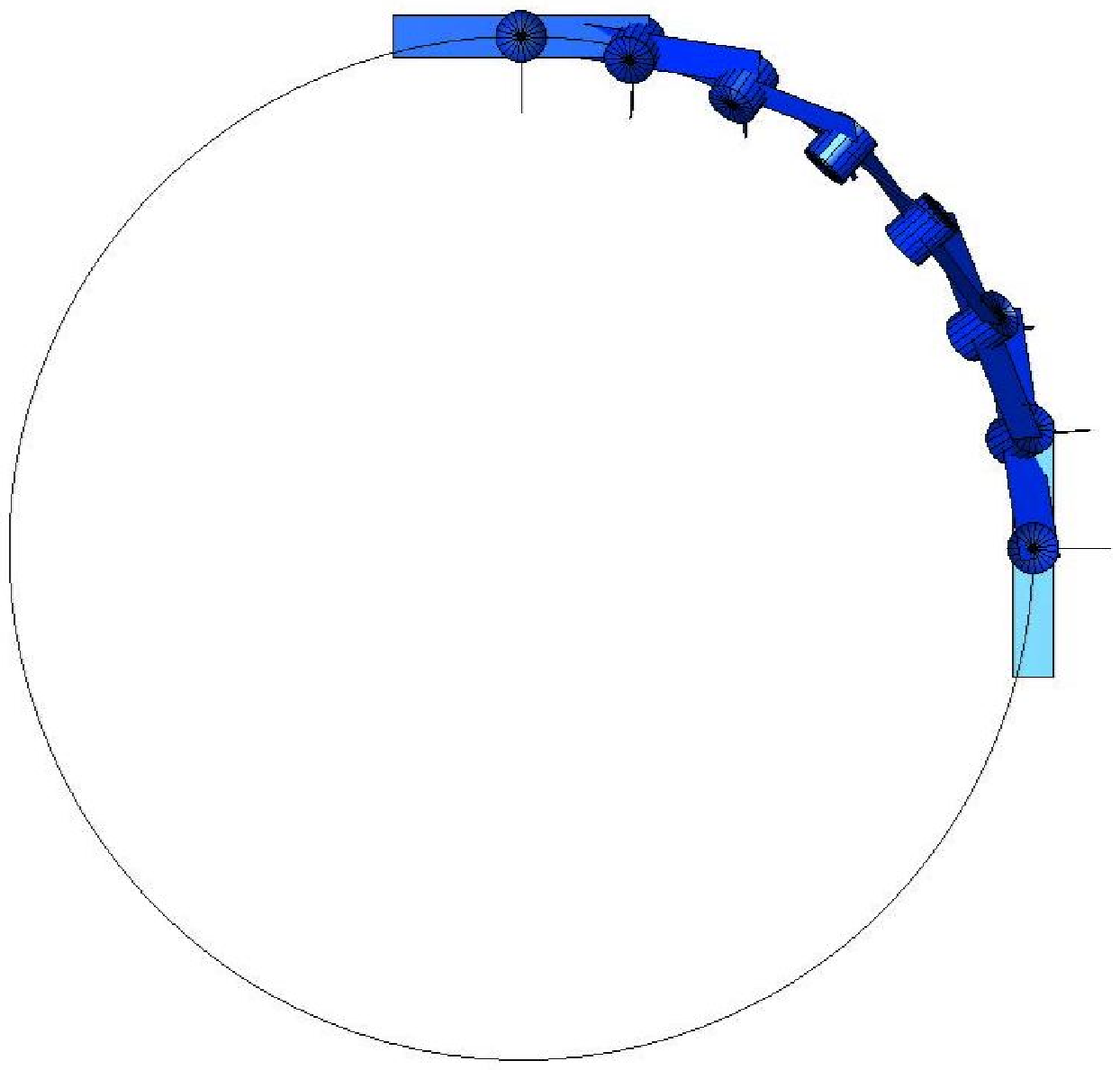}\label{fig:3di}
    \hfill
    \includegraphics[width=0.45\columnwidth]{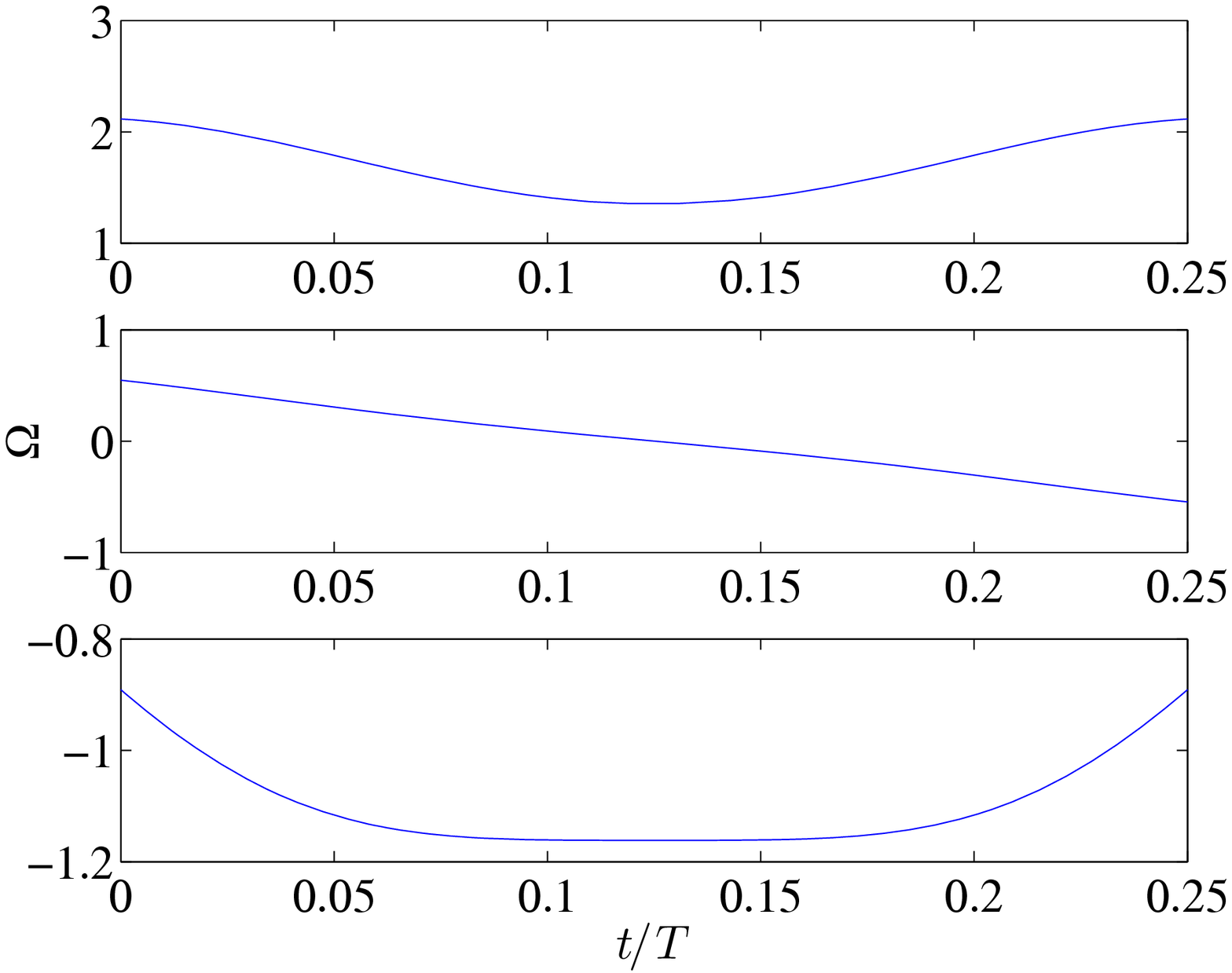}\label{fig:omegai}
    }
    \subfigure[Rotation about the LVLH axes $e_1$ and $e_2$]{
    \includegraphics[width=0.45\columnwidth]{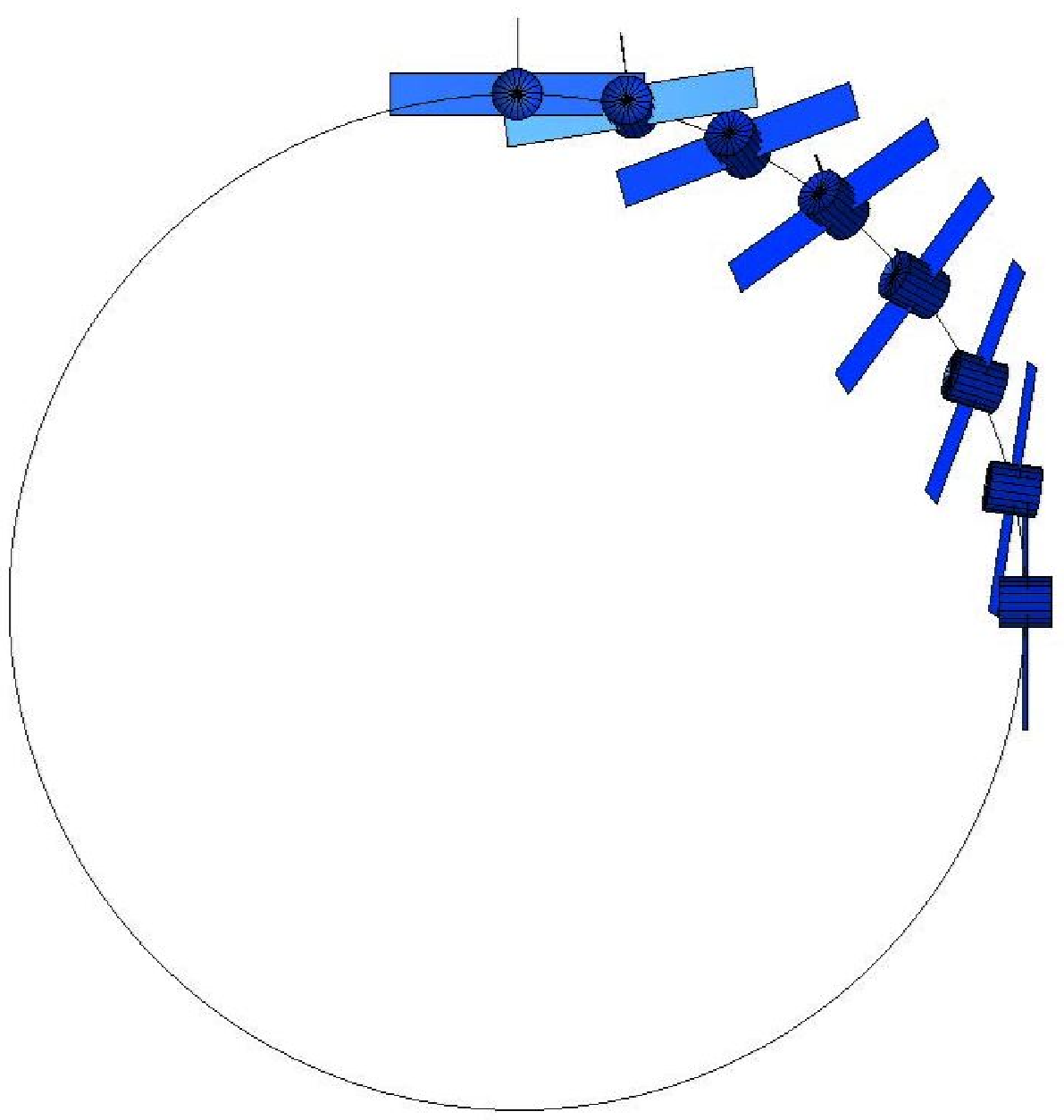}\label{fig:3dii}
    \hfill
    \includegraphics[width=0.45\columnwidth]{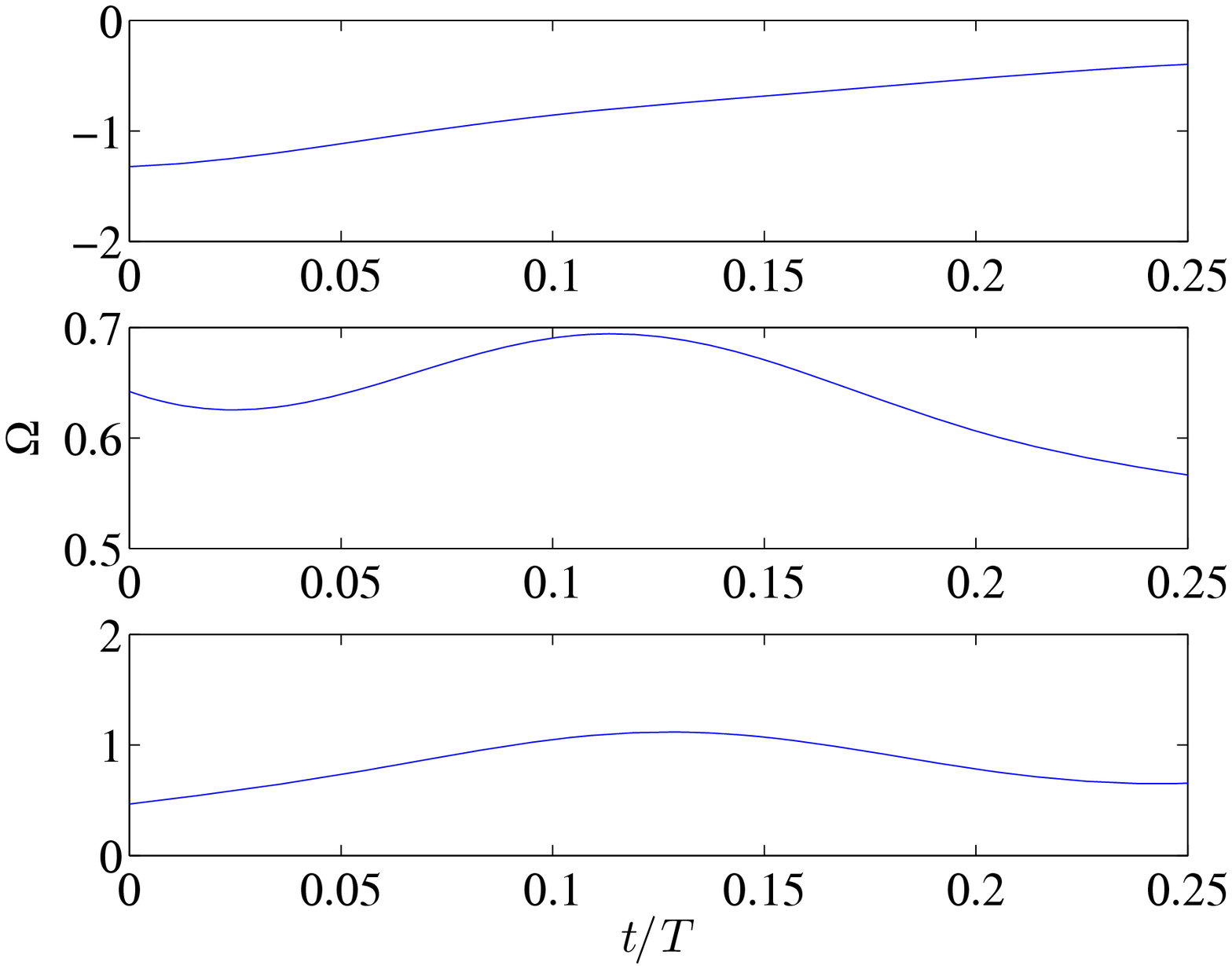}\label{fig:omegaii}
    }
    \subfigure[Rotation about the LVLH axes $e_2$ and $e_3$]{
    \includegraphics[width=0.45\columnwidth]{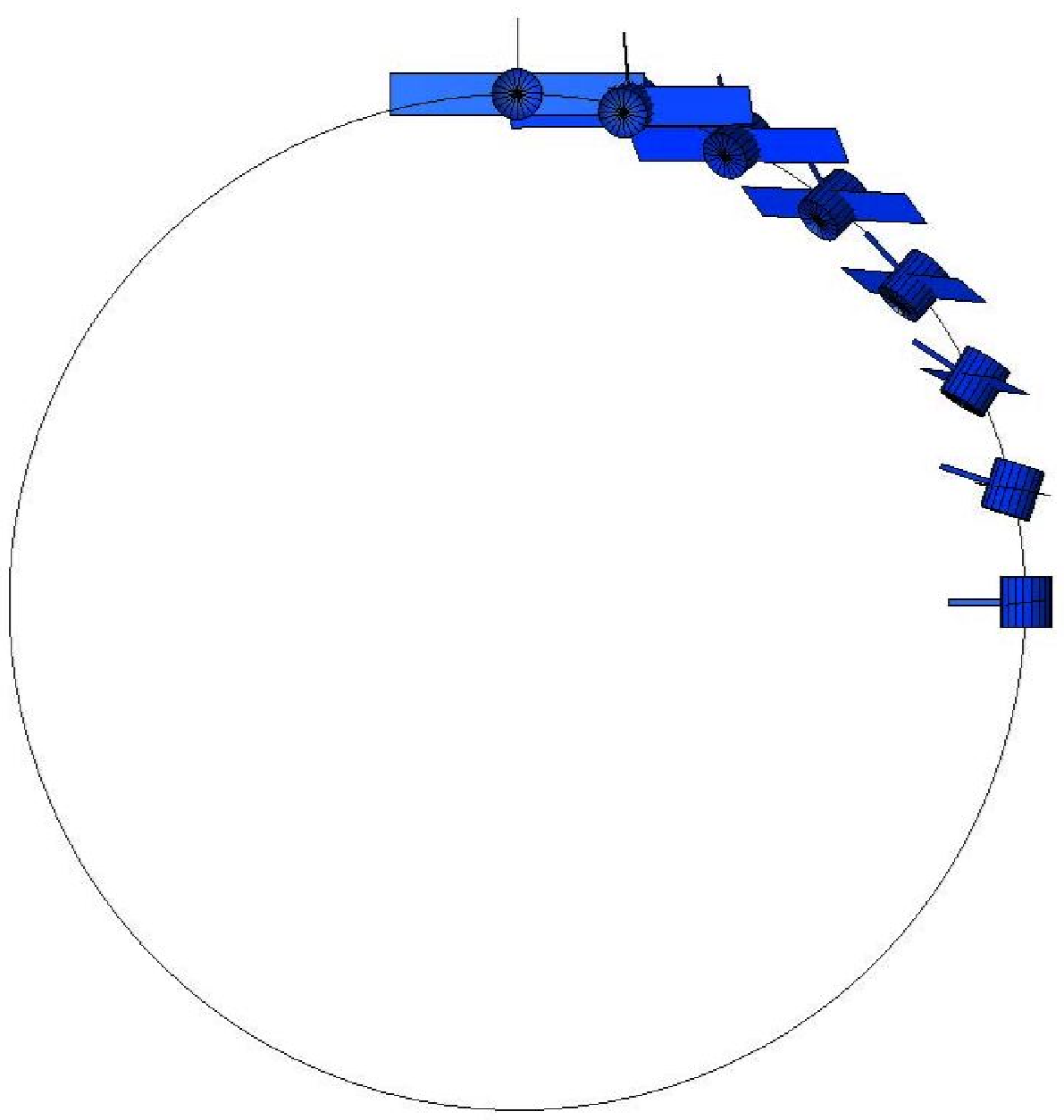}\label{fig:3diii}
    \hfill
    \includegraphics[width=0.45\columnwidth]{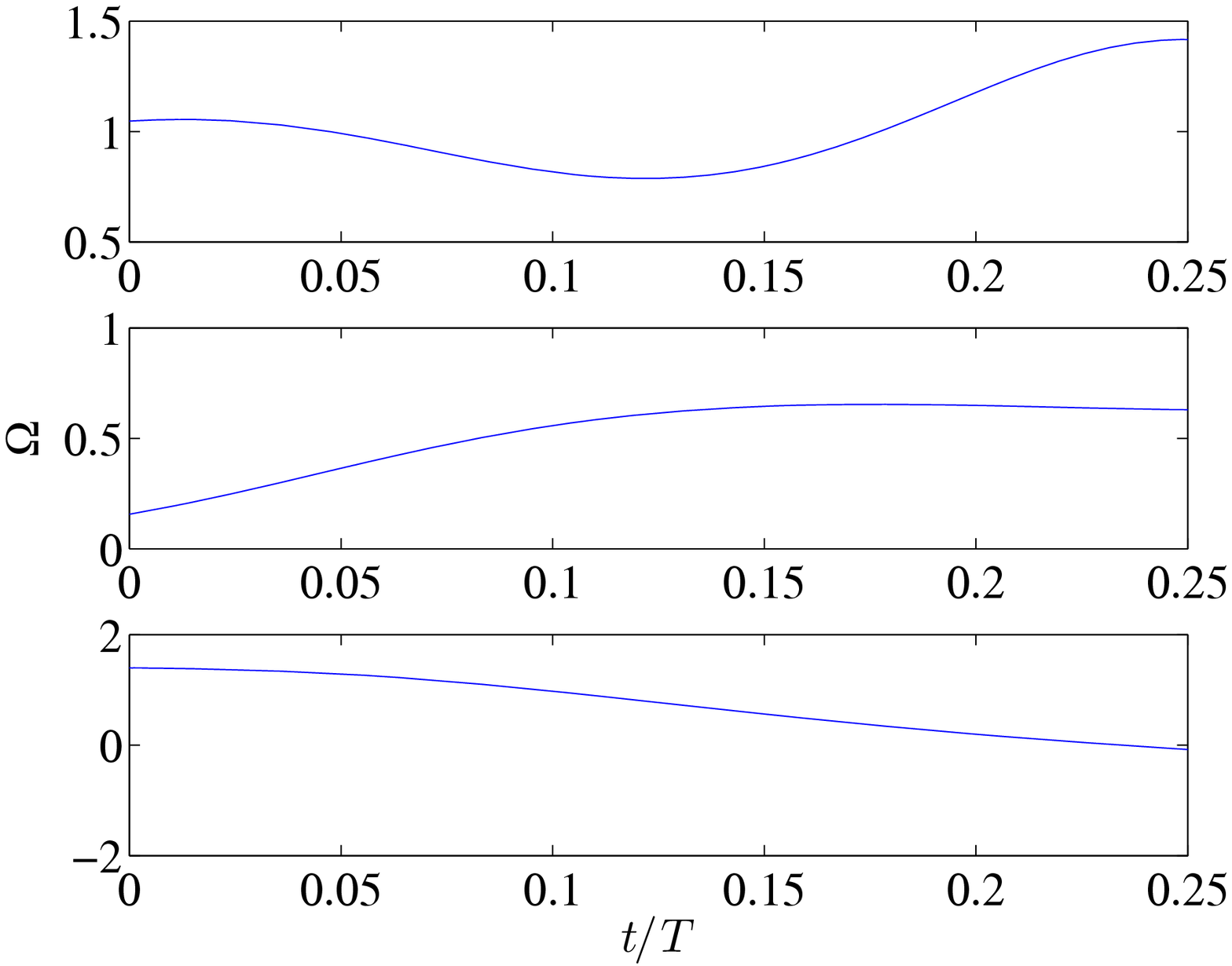}\label{fig:omegaiii}
    }
    \caption{Spacecraft attitude maneuvers}
    \label{fig:BVP}
\end{figure}

\reffig{BVP} shows the attitude maneuver of the spacecraft, and the
angular velocity response for each case. (Simple animations which show these spacecraft maneuvers can be found at \url{http://www.umich.edu/~tylee}.)

\section{Optimal spacecraft attitude maneuvers}\label{sec:opt}
\subsection{Problem formulation}
An optimization problem is formulated as a rest-to-rest maneuver of
an axially-symmetric spacecraft from a given initial attitude to a
given terminal reduced attitude for a fixed maneuver time. The
initial attitude is expressed by a rotation matrix with respect to
the LVLH frame, namely $R^{bl}_{0}$. The terminal desired attitude
is given by the reduced attitude, $\Lambda_{N_d}=R^{bl}_N
e_3\in\mathbb{S}^2$. This reduced attitude represents the direction
of the spacecraft axis of symmetry $e_3$ in the LVLH frame. Two
impulsive control moments are applied at the initial time and the
terminal time, and the maneuver of the spacecraft between the
initial time and the terminal time is uncontrolled.

In the problem studied in section \ref{sec:tpbvp}, the initial
parameter $\Pi_0$ is exactly prescribed by the constraint
$R^{bl}_{N_d}=R^{bl}_N$, and the discrete dynamics.   In this
section, we relax the terminal constraint by only specifying it up
to a rotation about the axis of symmetry of the spacecraft. Then, we can formulate an optimal attitude
maneuver problem.

The performance index is the sum of the magnitudes of the initial
impulse and the terminal impulse. Equivalently, one can minimize the
change in the initial angular momentum and the change in the terminal angular momentum. Since the initial attitude and the terminal time are fixed,
$\Pi_N$ and $\Lambda_N$ can be considered as functions of $\Pi_0$
through the discrete equations of motion. The
optimization problem is equivalent to
\begin{gather*}
\text{given : } R_0^{bl}, \Lambda_{N_d}, N,\\\
\begin{aligned}
\min_{\Pi_0} \mathcal{J} &= \norm{\Pi_0 - \omega_0 R_0^T  J e_2}
+\norm{\omega_0R_N^T  J e_2-\Pi_N},\\
&=\norm{H_0}+\norm{H_N},
\end{aligned}\\
\text{such that }
\mathcal{C}=\norm{\Lambda_{N}-\Lambda_{N_d}}^2=0,\\
\text{subject to \refeqn{updatePi}--\refeqn{updateRbl}}.
\end{gather*}

\addtolength{\textheight}{-0.4cm}

\subsection{Computational approach}
This problem is optimized by the Sequential Quadratic Programming
(SQP) method using analytical expressions for the sensitivity derivatives of the performance
index and of the constraint equation.

The variation of the performance index is
\begin{align*}
    \delta \mathcal{J} 
    &= \frac{H_0^T}{\norm{H_0}}\delta\Pi_0+\frac{H_N^T}{\norm{H_N}}\braces{-\omega_0S(\zeta_N)R_N^T  Je_2-\delta
    \Pi_N}.
\end{align*}
Since the initial attitude is given and fixed, the perturbation of
the initial attitude $\zeta_0$ is zero. Therefore
$\delta\Pi_N=\Phi_{22}\delta\Pi_0$ and
$\zeta_N=\Phi_{12}\delta\Pi_0$ from \refeqn{solsl}. Then,
$\delta\mathcal{J}$ is given by
\begin{align}\label{eqn:delJ}
    \delta \mathcal{J} =
    \bracket{\frac{H_0^T}{\norm{H_0}}+\frac{H_N^T}{\norm{H_N}}\braces{\omega_0S(R_N^TJe_2)\Phi_{12}-\Phi_{22}}}\delta\Pi_0.
\end{align}

Since $\Lambda_{N_d}$ is fixed and $\Lambda_N\in\S^2$, the variation of the constraint can be written as
\begin{align*}
    \delta \mathcal{C}
    & = -2\Lambda_{N_d}^T\delta\Lambda_N = 2\Lambda_{N_d}^T R_N^{bl}S(e_3)\zeta_N,
\end{align*}
where $\zeta_N=\Phi_{12}\delta\Pi_0$ from \refeqn{solsl}. Thus, $\delta \mathcal{C}$ is
\begin{align}\label{eqn:delg}
    \delta \mathcal{C} = \bracket{2\Lambda_d^T R_N^{bl} S(e_3)\Phi_{12}}\delta\Pi_0.
\end{align}
Equations \refeqn{delJ} and \refeqn{delg} are analytical expressions
for the sensitivity derivatives of the performance index and the
constraint.

\subsection{Numerical example}

The mass, length and time dimensions are normalized. The moment of
inertia of the spacecraft is chosen as
${J}=\mathrm{diag}\bracket{3,\,3,\,2}$, so that $e_3$ is the axis of
symmetry of the spacecraft.

The desired maneuver is to rotate the axis of symmetry from the
radial direction to the normal to the orbital plane during a quarter
orbit. The boundary conditions are given by
\begin{align*}
R_0^{bl} = \mathrm{diag}\bracket{1,\,-1,\,-1},\;
\Lambda_{N_d} = \bracket{0,\,-1,\,0}^T.
\end{align*}

\begin{figure}
    \subfigure{
    \includegraphics[width=0.45\columnwidth]{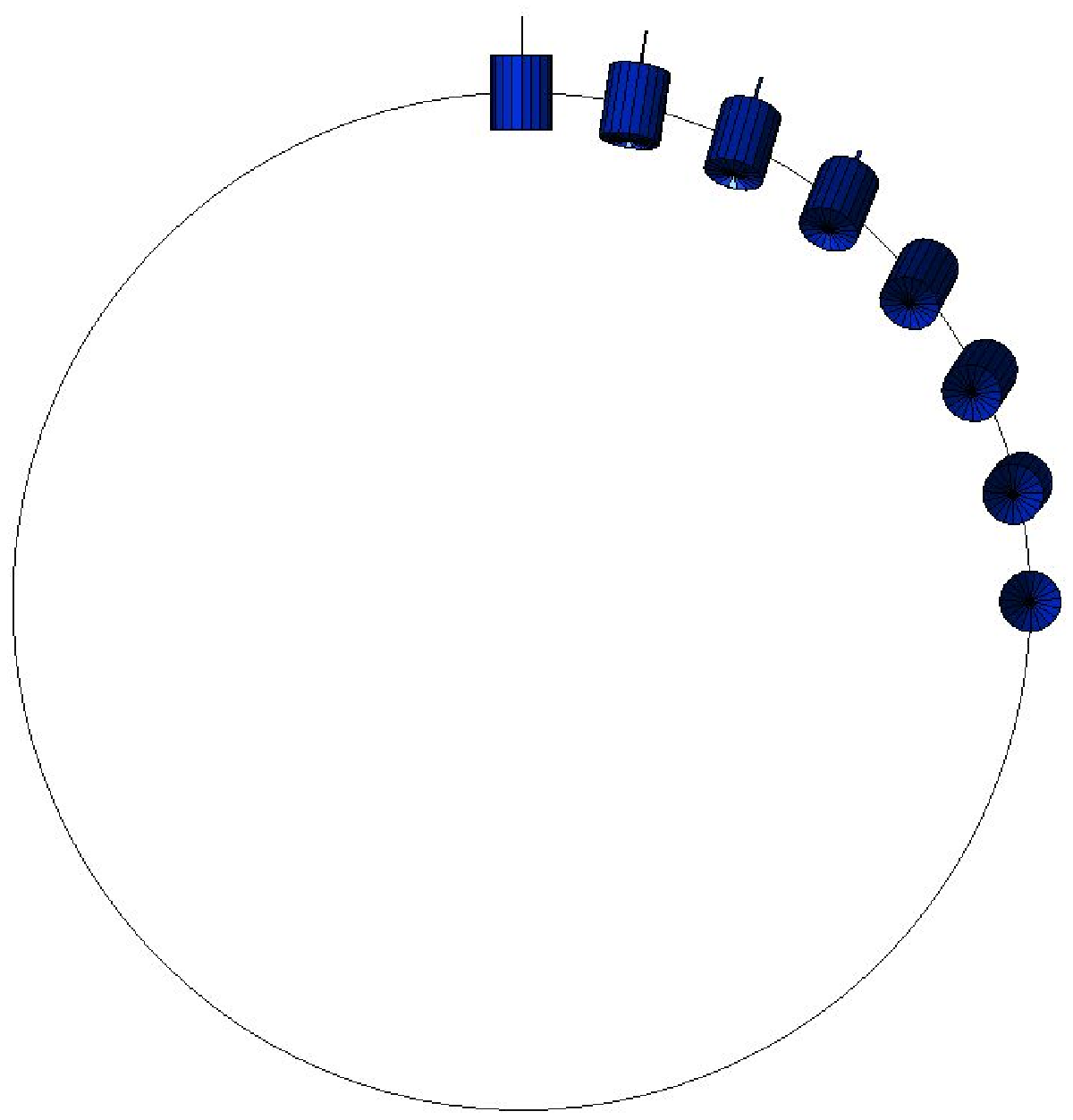}\label{fig:3dopt}
    \hfill
    \includegraphics[width=0.45\columnwidth]{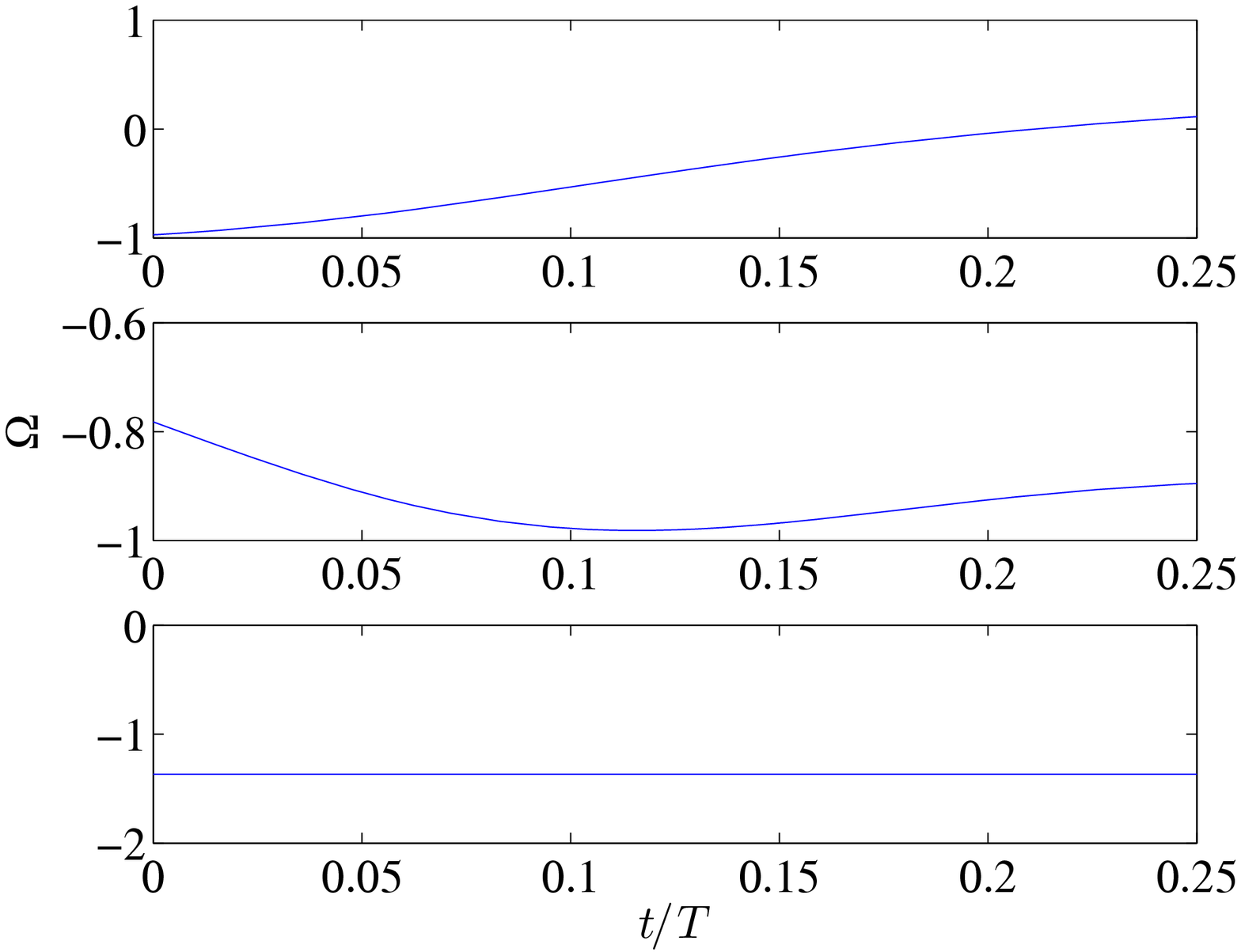}\label{fig:omegaopt}
    }
\caption{Optimal spacecraft attitude maneuver}
    \label{fig:opt}
\end{figure}

We use MATLAB's \verb"fmincon" function as an optimization tool. The
sensitivity derivatives of the performance index and the constraint
are provided by \refeqn{delJ} and \refeqn{delg}. The initial guess
of the initial angular momentum is chosen as $ {\Pi}_0^{(0)} =
{J}\bracket{1,\, 1,\, 0}^T$. The optimized performance index and the
violation of constraints are ${\mathcal{J}}=6.771$,
$\mathcal{C}=4.80\times 10^{-14}$. The corresponding angular momenta
and the terminal attitude are
\begin{gather*}
{\Pi}_0 = \bracket{-2.915,\,-2.347,\,-2.734}^T,\\
-{\Pi}_N = \bracket{-0.343,\,2.686,\,2.734}^T,\\
    R^{bl}_N =
    \begin{bmatrix}
   0.633&   0.733&    0.000\\
   0.000&    0.000&    -1.000\\
   -0.773&    0.633&   0.0000
   \end{bmatrix},
\end{gather*}
so that
$\Lambda_N=R_N^{bl}e_3=\bracket{0,\,0,\,-1}^T=\Lambda_{N_d}$.
\reffig{opt} shows the optimal maneuver of the spacecraft.

\section{Conclusion}
A global model for a rigid spacecraft in a circular orbit about a
large central body is presented. This model includes gravity
gradient effects that arise from the non-uniform gravity field.

The sensitivity derivatives for attitude dynamics of a rigid body
are derived while satisfying the global geometry of the problem.
Accurate computational approaches for solving a nonlinear boundary
value problem and the minimal impulse optimal control problem for
spacecraft attitude maneuvers are studied using sensitivity
derivatives.

The attitude dynamics are represented by a rotation matrix in the
Lie group $\SO$, and it is updated by Lie group variational
integrators that preserve the structure of $\SO$ as well as other geometric invariants of motion. The sensitivity derivatives are expressed in
terms of the Lie algebra $\so$. This approach completely avoids the
singularities and ambiguities associated with Euler angles or
quaternions, and it leads to a geometrically exact and numerically
efficient method for rigid body attitude dynamics problems.

Although the development in this paper includes a gravity gradient
moment and the rotation of the LVLH frame, the results presented
reduce to the case of a free rigid body if $\omega_0=0$. That is,
the computational approach suggested applies directly to attitude
maneuvers of the free rigid body.

\bibliography{ACC06}
\bibliographystyle{IEEEtran}

\end{document}